\documentclass[11pt]{article}

\usepackage{mathrsfs}
\usepackage{amssymb}
\usepackage{amscd}
\usepackage{amsfonts}

\font\tenmsb=msbm10
\font\sevenmsb=msbm7
\font\fivemsb=msbm5

\catcode`\@=11
\ifx\amstexloaded@\relax\catcode`\@=\active
\endinput\else\let\amstexloaded@\relax\fi
\def\spaces@{\space\space\space\space\space}
\def\spaces@@{\spaces@\spaces@\spaces@\spaces@\spaces@}
\def\space@.{\futurelet\space@\relax}
\space@.
\def\Err@#1{\errhelp\defaulthelp@\errmessage{AmS-TeX error: #1}}
\def\relaxnext@{\let\next\relax}
\def\accentfam@{7}
\def\noaccents@{\def\accentfam@{0}}
\def\Cal{\relaxnext@\ifmmode\let\next\Cal@\else
\def\next{\Err@{Use \string\Cal\space only in math mode}}\fi\next}
\def\Cal@#1{{\Cal@@{#1}}}
\def\Cal@@#1{\noaccents@\fam\tw@#1}
\def\Bbb{\relaxnext@\ifmmode\let\next\Bbb@\else
\def\next{\Err@{Use \string\Bbb\space only in math mode}}\fi\next}
\def\Bbb@#1{{\Bbb@@{#1}}}
\def\Bbb@@#1{\noaccents@\fam\msbfam#1}
\newfam\msbfam
\textfont\msbfam=\tenmsb
\scriptfont\msbfam=\sevenmsb
\scriptscriptfont\msbfam=\fivemsb

\def\NN{{\Bbb N}}
\def\ZZ{{\Bbb Z}}
\def\RR{{\Bbb R}}

\def\beq{\begin{equation}}
\def\eeq{\end{equation}}

\def\skipaline{\removelastskip\vskip12pt plus 1pt minus 1pt}

\def\Proof{\removelastskip\skipaline
\noindent \it Proof. \rm}

\newtheorem{Theorem}{Theorem}
\newtheorem{Lemma}{Lemma}[section]

\newtheorem{Remark}{Remark}[section]

\@addtoreset{equation}{section}

\begin{document}

\title{Boundedness of solutions for a class of impact oscillators with time-denpendent polynomial potentials%%\footnote{}
\author{Daxiong \ Piao\footnote{Supported by the NSF of Shangdong Province(No.ZR2012AM018),  E-mail: dxpiao@ouc.edu.cn
}, \,\,\,\,\,Xiang \ Sun\footnote{E-mail: sxltlg2312@126.com}
 \\
        School of Mathematical Sciences,\ Ocean University of China
\\
        \ \ Qingdao 266100,\ \ P.R. China
}
\date{}}

\maketitle

\begin{abstract}In this paper, we consider the boundedness of solutions for a class of impact oscillators with time dependent polynomial potentials,\\
%$$
%\left\{\begin{array}{11}
%\ddot{x}+x^{2n+1}+\sum_{i=0}^{2n}p_{i}(t)x^{i}(t)& \quad {\rm for}\quad x(t)> 0,\\
%x(t)\geq 0,\\
%\dot{x}(t_{0}^{+}=-\dot{x}(t_{0}^{-}& \quad {\rm if}\quad x(t_{0})=0,
%\end{array}\right.
%$$
$$
 \left\{\begin{array}{ll}
 \displaystyle
 \ddot{x}+x^{2n+1}+\sum_{i=0}^{2n}p_{i}(t)x^{i}=0,& \quad {\rm for}\quad x(t)> 0,\\
 x(t)\geq 0,&\\
 \dot{x}(t_{0}^{+})=-\dot{x}(t_{0}^{-}),& \quad {\rm if}\quad x(t_{0})=0,
 \end{array}\right.
 $$
where $n\in\NN^{+}$, $p_{i}(t+1)=p_{i}(t)$ and $p_{i}(t)\in C^5(\RR/\ZZ).$
%AMS(1991) Subject Classification:  34C99, 58F10.
\end{abstract}

\vskip 1.0cm
\section{Introduction and main result}

\ \ \ \ \ \ In Ref. \cite{[DZ]}, Dieckerhoff and Zehnder proved
the boundedness of solutions for a time-dependent nonlinear
differential equation: \beq\label{DZ}
\ddot{x}+x^{2n+1}+\sum_{i=0}^{2n}p_{i}(t)x^{i}=0, \quad n\geq 1,
\eeq with $p_{i}(t+1)=p_{i}(t)$ and $p_{i}(t)\in {\cal C}^\infty$
and asked whether or not the boundedness phenomenon is related to
the smoothness of $p_{i}(t).$ Then Laederich and Levi \cite{[LL]}
relaxed the smoothness requirement of $p_{i}(t)$ to ${\cal
C}^{5+\varepsilon}$ with $\varepsilon>0,$ Yuan \cite{[Y1],[Y2]}
relaxed the smoothness requirement of $p_{i}(t).$ to ${\cal C}^2,$
etc.

Motivated by \cite {[DZ],[LL],[Y1],[Y2]}, we are going to consider the problem of the boundedness for Eq.(\ref{DZ}) with impacts:
\beq\label{imp}
 \left\{\begin{array}{ll}
 \displaystyle
 \ddot{x}+x^{2n+1}+\sum_{i=0}^{2n}p_{i}(t)x^{i}=0,& \quad {\rm for}\quad x(t)> 0,\\
 x(t)\geq 0,&\\
 \dot{x}(t_{0}^{+})=-\dot{x}(t_{0}^{-}),& \quad {\rm if}\quad x(t_{0})=0,
 \end{array}\right.
\eeq
where $n\in\NN^{+}$, $p_{i}(t+1)=p_{i}(t)$ and $p_{i}(t)\in C^5(\RR/\ZZ).$

The system is included in the system of impact oscillators given by
\beq\label{osc}
 \left\{\begin{array}{ll}
 \ddot{x}+V_{x}(x,t)=0,& \quad {\rm for}\quad x(t)> 0,\\
 x(t)\geq 0,&\\
 \dot{x}(t_{0}^{+})=-\dot{x}(t_{0}^{-}),& \quad {\rm if}\quad x(t_{0})=0,
 \end{array}\right.
\eeq
which serve as models of dynamical systems with discontinuities \cite{[Kun]}. From the viewpoint of mechanics,
this equation models the motion of particle attached to a nonlinear spring and bouncing elastically against the
fixed barrier. The system of this form also relate to the research of the Fermi accelerator \cite{[LC]},
dual billiards \cite{[Boy]} and celestial mechanics \cite{[CL]}.

The nonsmoothness caused by the impact limit the applications of many powerful mathematical tools.
However there are also many interesting papers on the impact oscillators, see \cite{[QT1]}-\cite{[WRQ]} and their references.

To overcome the problem of nonsmoothness, Qian and Torres \cite{[QT1]} and Qian and Sun \cite{[QX]}
used success mapping in their papers and proved the existence of invariant tori by a variant version
of Moser's small twist theorem by Ortega \cite{[Or]}. In Zharnitsky \cite{[Zha]}, Wang \cite{[WW],[WLQ]},
they exchanged the position of angle and time variables of the Hamiltonian function which only ${\cal C}^0$
in angle variable, and finally obtained the boundedness of solutions by Moser's small twist theorem.

In \cite{[WW]}, Z.Wang and Y.Wang studied the following impact
oscillator: \beq\label{WW}
 \left\{\begin{array}{ll}
 \ddot{x}+x^{2n+1}=p(t),& \quad {\rm for}\quad x(t)> 0,\\
 x(t)\geq 0,&\\
 \dot{x}(t_{0}^{+})=-\dot{x}(t_{0}^{-}),& \quad {\rm if}\quad x(t_{0})=0,
 \end{array}\right.
\eeq
where $n\in\NN^+,$ $p(t)$ is a ${\cal C}^5$ periodic function with 1, and the term $x^{2n+1}$ models a
hard spring. They proved the Lagrangian stability for system (\ref{WW}) by exchanging the roles of angle
and time variables and using Moser's small twist theorem.

In this paper we will extend (\ref{WW}) to the case of time-periodic polynomial potentials with constant
leading coefficient as (\ref{imp}) and use similar method to prove the boundedness of solutions for (\ref{imp}).
More exactly, we obtain the following conclusion:
%\begin{theorem} Every solution of (1.2) is bounded, i.e. $x(t)$ exits for $t\in\RR$ and
%$sup(|x(t)|+|\dot{x}(t)|)<+\infty.$
%\end{theorem}
\begin{Theorem} Every solution of (1.2) is bounded, i.e., $x(t)$ exits for $t\in\RR$ and\\
$$\sup(|x(t)|+|\dot{x}(t)|)<+\infty.$$%\end{equation}
\end{Theorem}

The idea of proving the boundedness of solutions of (\ref{imp}) is as follows. First of all,
we describe (\ref{imp}) by a Hamiltonian function $H_3(\rho,\phi,t)$ (see (\ref{hami3}))
in action-angle variables defined on the whole space $\RR^+\times S^1\times S^1$ by means
of transformation theory. Due to the existence of impact, $H_3(\rho,\phi,t)$ is only
continuous in $\phi$, but is ${\cal C}^5$ smooth in $\rho$ and $t$. Exchanging the roles
of the variables $\phi$ and $t$ outside of a large disk $D_{r}=\{(\rho,\phi),\rho<r\}$ in $(\rho,\phi)$-plane, (2.9)
is transformed into a perturbation of an integrable Hamiltonian system $H_4(I,\theta,\tau)$ (see (3.3))
which is sufficiently smooth in $I$ and $\theta$. The Poincar\'{e} mapping with respect to the new time $\tau$
is closed to a so-called twist mapping in $(\RR^+\times S^1)\setminus D_r$, and satisfies Moser's invariant curve
theorem after a scaling transformation. Then Moser's theorem guarantees the existence of arbitrarily large invariant
curves diffeomorphic to $\upsilon=const.$ over  $\upsilon=0$ in  $(\upsilon,\theta)$-plane. Go back to the equivalent
system (2.1), every such curve is a base of a time-periodic and flow-invariant cylinder in the extended phase
space $(x,y,t)\in \RR^+\times \RR\times \RR$, which confines the solutions in the interior and which leads to a
bound of these solutions if the uniqueness of initial value problem holds. Note that, from the proof below,
system (2.1) is equivalent to a smooth system (2.9) in $(\RR^+\times S^1)\setminus D_r$ which has uniqueness.
Hence system (2.1) has uniqueness.
\begin{Remark} The time-dependence of the stiffness coefficient $p_i(t)$ of the spring can be produced by periodic
changes of the temperature or other physical variables. Following
Laederich and Levi \cite{[LL]} and Z.Wang and Y. Wang \cite{[WW]},
we assume $p_i(t)\in {\cal C}^5$.\end{Remark}
\begin{Remark} Compared with (1.4), the potentials of (1.1) become more complicated which will make more difficulties in the proof of Theorem 1.\end{Remark}
\begin{Remark} The main difficulty in this paper is also the nonsmoothness caused by the impact. Similar to \cite{[WW]}, we will also exchange the roles of angle and time variables of the Hamiltonian function $H_3(\rho,\phi,t)$ (see below), and then obtain the boundedenes of solutions by Moser's small twist theorem.\end{Remark}

The rest of this paper is organized as follows. In section 2, we will make some action-angle transformations to transform the system into an equivalent Hamiltonian $H_4(I,\theta,\tau)$ (see below) which is defined in the whole plane $\RR^+\times S^1\times S^1$ and ${\cal C}^5$ in $I,\theta$, but only ${\cal C}^0$ in $\tau$. In section 3, we do more canonical transformations such that the Poincar\'{e} mapping of the new system is closed to the twist mapping and then use Moser's theorem to complete the proof of Theorem 1. Finally, we will give the proofs of some technical lemmas in section 4 and draw a conclusion in section 5.

Throughout this paper, we denote $A<.B$ if $A<C_1B$, and $A>.B$ if $A>C_2B$ , where $C_1,C_2$ are positive constants.

\vskip 2mm
\section{Action-angle variables}

\ \ \ \ \ \ In this section, we will make some action-angle transformations to transform the system into an equivalent Hamiltonian $H_4(I,\theta,\tau)$ (see below) which is defined in the whole plane $\RR^+\times S^1\times S^1$ and ${\cal C}^5$ in $I,\theta$, but only ${\cal C}^0$ in $\tau$.

Without impact, Eq.(\ref{imp}) is just Eq.(\ref{DZ}) which is of second order and equivalent to the one order system
\beq\label{sys1}
 \left\{\begin{array}{ll}
 \dot{x}=y,\\
 \displaystyle
 \dot{y}=-x^{2n+1}-\sum_{i=0}^{2n}p_i(t)x^i.
 \end{array}\right.
\eeq
(2.1) is a system defined in the whole phase plane $XOY$, which has a Hamiltonian function
\beq\label{hami}
H(x,y,t)=\frac{1}{2}y^2+\frac{1}{2n+2}x^{2n+2}+\sum_{i=0}^{2n}\frac{p_i(t)x^{i+1}}{i+1}
\eeq
with the symplectic form $dx\wedge dy.$

In order to make an action-angle transformation, we consider the system
\beq\label{sys2}
 \left\{\begin{array}{ll}
 \dot{x}=y,\\
 \dot{y}=-x^{2n+1}.
 \end{array}\right.
\eeq
with the Hamiltonian $H_0(x,y,t)=\frac{1}{2}y^2+\frac{1}{2n+2}x^{2n+2}.$

Clearly, $H_0>0$ on $\RR^2$ except at the only equilibrium point $(x,y)=(0,0)$ where $H_0=0$.

Let $(C(t),S(t))$ is the solution of Eq.(2.3) satisfying the initial condition:$(C(0),S(0))=(1,0)$. Let $T_0>0$ be its minimal period.

From Eq.(2.3), we can find that $C(t)$ and $S(t)$ satisfy:
$$
\begin{array}{ll}
{\rm (1)}&C(t),\ S(t)\in {\cal C}^\infty(\RR), \ and\ C(t+T_0)=C(t),\ S(t+T_0)=S(t)\ with\ C(0)=1,\ S(0)=0;\\
{\rm (2)}&\dot{C}(t)=S(t)\ and\ \dot{S}(t)=-C(t)^{2n+1};\\
{\rm (3)}&(n+1)S(t)^2+C(t)^{2n+1}=1;\\
{\rm (4)}&C(-t)=C(t)\ and\ S(-t)=-S(t).
\end{array}
$$

Now we can make coordinate transformation by the mapping $\Psi_1:\RR^+\times S^1\rightarrow \RR^2/\{0\},$ where $(x,y)=\Psi_1(\lambda,\vartheta)$ defined by the formula
$$
x=(a\lambda)^\alpha C(\vartheta T_0),\ \ \ y=(a\lambda)^\alpha S(\vartheta T_0),
$$
where $\alpha=\frac{1}{n+2},\ \beta=1-\alpha,\ a=\frac{T_0}{\alpha}$ are constants. We can claim that $\Psi_1$ is a symplectic diffeomorphism from $\RR^+\times S^1$ onto $\RR^2/\{0\}$. Indeed, the Jacobian determinant of $\Psi_1$ is 1, so $\Psi_1$ is measure preserving. Moreover, since $(C(t),S(t))$ is a solution of Eq.(\ref{sys2}) and has $T_0$ as its minimal period, which concludes that $\Psi_1$ is one to one and onto. This proves the claim.

Under $\Psi_1$, Hamiltonian $H$ is transformed into
\beq\label{hami1}
H_1(\lambda,\vartheta,t)=\frac{a^{2\beta}}{2n+2}\lambda^{2\beta}+\sum_{i=0}^{2n}\frac{p_i(t)(a\lambda)^{\alpha(i+1)}[C(\vartheta T_0)]^{i+1}}{i+1},
\eeq
which is ${\cal C}^5$ in $\lambda$, $t$ and ${\cal C}^0$ in $\vartheta.$

Considering the impact case, the phase space is only a half plane $(x\geq 0)$ of the original phase plane $XOY$. Since impact, the smoothness with variable $\vartheta$ of $H_1$ is disappeared. We shall study another Hamiltonian $H_3(\rho,\phi,t)$ as defined below.  Let $C(\vartheta T_0)$, we will find two solutions $\vartheta=\vartheta_1, \vartheta_2$ in interval $[0,1)$. By the symmetry of $H_0$, we can easily find $\vartheta_1=0, \vartheta_2=\frac{1}{2}.$

Under $\Psi_1$, $(\lambda,0(mod\ 1))$ is mapped into the positive $Y$-axis, $(\lambda,\frac{1}{2}(mod\ 1))$ is mapped into the negative $Y$-axis.

Define a sympletic transformation $\Psi_2:(\phi,\rho)\mapsto(\vartheta,\lambda)$, given by $\vartheta=\frac{1}{2},\ \ \ \lambda=2\rho.$ And $H_1$ is transformed into
\beq\label{hami2}
H_2(\rho,\phi,t)=d\rho^{2\beta}+\sum_{i=0}^{2n}\frac{p_i(t)(2a\rho)^{\alpha(i+1)}[C(\frac{1}{2}\phi T_0)]^{i+1}}{i+1},
\eeq
where $d=\frac{(2a)^{2\beta}}{2n+2}.$

To get a system which is equivalent to (1.2), we could define a new Hamiltonian function by
\beq\label{hami3}
H_3(\rho,\phi,t)=H_2(\rho,\phi-[\phi],t)=d\rho^{2\beta}+\sum_{i=0}^{2n}\frac{p_i(t)(2a\rho)^{\alpha(i+1)}\{C(\frac{1}{2}(\phi-[\phi]) T_0)\}^{i+1}}{i+1},
\eeq
where $[\phi]$ denotes the largest integer less than or equal to $\phi$. And the corresponding Hamiltonian system is
\beq\label{sys3}
 \left\{\begin{array}{ll}
 \dot{\rho}=\frac{\partial H_3}{\partial \phi},\quad {\rm when}\quad \phi(t)\in (k,k+1), k\in\ZZ,\\
 \dot{\phi}=-\frac{\partial H_3}{\partial \rho},\quad {\rm when}\quad \phi(t)\in (k,k+1), k\in\ZZ,\\
 \displaystyle
 \phi(t_0)=k,\ \rho(t_0)=\lim_{t\rightarrow t_0}\rho(t),\quad {\rm when}\quad \lim_{t\rightarrow t_0}\phi(t)=k, k\in\ZZ.
 \end{array}\right.
\eeq

Obviously, $H_3$ is periodic in $\phi$ with 1, and ${\cal C}^5$ in $\rho, \phi$ when $\phi\notin\ZZ$, but only continuous in $\phi$ when $\phi\in\ZZ$. In fact, $H_3$ is right continuous when $\phi\in\ZZ.$

%\begin{Lemma} For every solution $(\rho(t),\phi(t))$ with $\rho(t)\neq 0$, $(x(t),\dot{x}(t))=\Psi_1\circ\Psi_2(\rho(t),\phi(t)-[\phi(t)])$ is a continuous solution of (1.1) with $(x(t)\neq 0$, and vice versa.\end{Lemma}
\begin{Remark} The systems mentioned above are equivalent, see \cite{[WW]} and \cite{[WLQ]}, we omit the proof.\end{Remark}

\vskip 2mm
\section{The proof of Theorem 1}
\ \ \ \ Now we are concerned with the Hamiltonian system (\ref{sys3}) with Hamiltonian function $H_3$ given by (\ref{hami3}). To cope with the nonsmoothness in $\phi$ and use Moser's small twist theorem to prove the Lagrangian stability, we will exchange the positions of variables $(\rho,\phi)$ and $(H_3,t)$ below. This trick has been used in \cite{[Lev],[Liu],[JPW]}.

From (\ref{hami3}), we have that
$$\lim_{\rho\rightarrow \infty}H_3=\infty$$
and
$$\frac{\partial H_3}{\partial \rho}=2d\beta\rho^{2\beta-1}+\sum_{i=0}^{2n}\frac{\alpha(i+1)p_i(t)(2a)^{\alpha(i+1)}\rho^{\alpha(i+1)-1}\{C(\frac{1}{2}(\phi-[\phi])T_0)\}^{i+1}}{i+1}>0,$$
if $\rho\geq 1.$ By the implicit theorem, we know that there is a function $R=R(H_3,t,\phi)$ such that
\beq\label{rid}
\rho(H_3,t,\phi)=d^{-\frac{1}{2\beta}}H_3^{\frac{1}{2\beta}}-R(H_3,t,\phi).
\eeq

Then we make another transformation $\Psi_3:(\rho,\phi,t)\mapsto(I,\theta,\tau)$ given by
\beq\label{map2}
I=H_3(\rho,\phi,t),\quad \theta=t,\quad \tau=\phi.
\eeq

This transformation also leads to a Hamiltonian system with new Hamiltonian function
\beq\label{hami4}
H_4(I,,\theta,\tau)=d^{-\frac{1}{2\beta}}I^{\frac{1}{2\beta}}-R(I,\theta,\tau),
\eeq
which is ${\cal C}^5$ in $\theta$, ${\cal C}^\infty$ in $I$ and ${\cal C}^0$ in $\tau.$

And the system is given by
\beq\label{sys4}
 \left\{\begin{array}{ll}
 \dot{I}=-\frac{\partial H_4}{\partial \theta},&\quad \tau\notin\ZZ,\\
 \dot{\theta}=\frac{\partial H_4}{\partial I},&\quad \tau\notin\ZZ,\\
 \displaystyle\theta(k)=\lim_{\tau\rightarrow k}\theta(\tau),\quad I(k)=\lim_{\tau\rightarrow k}I(\tau),&\quad k\in\ZZ.
 \end{array}\right.
\eeq

Then we have\\
\begin{Lemma} For $I$ large enough, the following estimate about $R$ hold ture
\beq\label{ineq1}
|D_I^jD_\theta^j R|<.I^{\frac{1}{2}-1},\quad {\rm for} \quad j+k\leq 5.
\eeq\end{Lemma}

To be convenient for readers, the proof is given in section 4.

In order to use the Moser's small twist theorem, we introduce a new variable $\upsilon$ and a small positive $\varepsilon$ by the formula $I=\frac{\upsilon}{\varepsilon}$, where $\upsilon\in [1,2]$, $\varepsilon>0$ and $\varepsilon\rightarrow 0$, when $I\rightarrow\infty$. Then the system (\ref{sys4}) is equivalent to the following system
\beq\label{sys5}
 \left\{\begin{array}{ll}
 \dot{\upsilon}=\varepsilon\frac{\partial R}{\partial \theta},&\quad \tau\notin\ZZ,\\
 \dot{\theta}=\frac{1}{2\beta}d^{-\frac{1}{2\beta}}\varepsilon^{1-\frac{1}{2\beta}}\upsilon^{\frac{1}{2\beta}-1}-\varepsilon\frac{\partial R}{\partial\upsilon},&\quad \tau\notin\ZZ,\\
 \displaystyle\theta(k)=\lim_{\tau\rightarrow k}\theta(\tau),\quad \upsilon(k)=\lim_{\tau\rightarrow k}\upsilon(\tau),&\quad k\in\ZZ.
 \end{array}\right.
\eeq

Integrate the system (\ref{sys5}) from $\tau=0$ to $\tau=1$, we obtain the Poincar\'{e} mapping $P$ of the form
\beq\label{Poinc}
\upsilon_1=\upsilon_0+f_1(\upsilon_0,\theta_0),\quad \theta_1=\theta_0+\frac{1}{2\beta}d^{-\frac{1}{2\beta}}\varepsilon^{1-\frac{1}{2\beta}}\upsilon_0^{\frac{1}{2\beta}-1}+f_2(\upsilon_0,\theta_0),
\eeq
where $f_1(\upsilon_0,\theta_0)=\int_0^1\varepsilon\frac{\partial R}{\partial\theta}d\tau,\quad f_2(\upsilon_0,\theta_0)=-\int_0^1\varepsilon\frac{\partial R}{\partial\upsilon}d\tau,$ which satisfy:
\begin{Lemma}
If $\varepsilon$ small enough, then
\beq\label{ineq2}
|D_{\upsilon_0}^jD_{\theta_0}^k f_1(\upsilon_0,\theta_0)|<.\varepsilon^\frac{1}{2},\quad {\rm for}\quad j+k\leq 4;
\eeq
and\\
\beq\label{ineq3}
|D_{\upsilon_0}^jD_{\theta_0}^k f_2(\upsilon_0,\theta_0)|<.\varepsilon^\frac{1}{2},\quad {\rm for}\quad j+k\leq 4.
\eeq
\end{Lemma}
\Proof Note $I=\upsilon/\varepsilon$, we have $D_{\upsilon}^jR(\upsilon/\varepsilon,\theta,\tau)=\varepsilon^{-j}D_{I}^{j}R(I,\theta,\tau),$ and $\varepsilon\rightarrow 0$, when $I\rightarrow\infty$. Thus, with (\ref{ineq1})-(\ref{Poinc}), it is readily to get the results, see \cite{[DZ]}. We omit it.

\begin{Lemma}
The mapping $P$ possesses the intersection property in the annulus $[1,2]\times S^1$, i.e., if $\Gamma$ is an embedded circle in $[1,2]\times S^1$ homotopic to a circle $\upsilon=const.$, then $P(\Gamma)\cap\Gamma\neq\emptyset.$
\end{Lemma}

The proof can be found in \cite{[DZ]}.\\
{\it Proof of Theorem 1.}\quad Until now, we have verified that the mapping $P$ satisfies all the conditions of Moser's small twist theorem \cite{[Mos],[Rus]}. Hence for any $\omega\in(\frac{1}{2\beta}d^{-\frac{1}{2\beta}}\varepsilon^{1-\frac{1}{2\beta}},\frac{1}{2\beta}d^{-\frac{1}{2\beta}}2^{1-\frac{1}{2\beta}})$  satisfying
\beq\label{ineq4}
|\omega-\frac{p}{q}|\geq.\varepsilon|q|^{-5/2},
\eeq
there is an invariant curve $\Gamma$ of $P$ surrounding $\upsilon_0=1$ if $\varepsilon$ is sufficiently small. This means that the Poincar\'{e} mapping of the system (\ref{sys4}) has an invariant curve $\widetilde{\Gamma}$ surrounding $I=\frac{1}{\varepsilon}$ and in the annulus $\frac{1}{\varepsilon}<I<\frac{2}{\varepsilon}.$ Then transform the system (\ref{sys4}) back to the system (\ref{sys3}), we can find that there exit arbitrary large invariant tori in $(\rho,\phi,t)$ space, and this implies that all the solutions are bounded. Finally, according to the equivalent of the system (\ref{DZ}) and (\ref{sys3}), we can imply the boundedness of all the solutions of Eq.(\ref{DZ}). The proof of Theorem 1 is completed.

\vskip 2mm
\section{The proof of Lemma 3.1}

\Proof From (\ref{rid}) and (\ref{map2}), we have
\beq\label{rid2}
\rho=d^{-\frac{1}{2\beta}}I^{\frac{1}{2\beta}}-R
\eeq
where $\alpha=\frac{1}{n+2},\ \beta=1-\alpha=\frac{n+1}{n+2},\ a=\frac{T_0}{\alpha}$ and $d=\frac{(2a)^{2\beta}}{2n+2}.$

When $j+k=0,$ substitute $\rho$ in (\ref{hami3}), we have
\beq\label{eq1}
I=d(d^{-\frac{1}{2\beta}}I^{\frac{1}{2\beta}}-R)^{2\beta}+\sum_{i=0}^{2n}\frac{p_i(\theta)(2a)^{\alpha(i+1)}(d^{-\frac{1}{2\beta}}I^{\frac{1}{2\beta}}-R)^{\alpha(i+1)}\{C(\frac{1}{2}(\tau-[\tau])T_0)\}^{i+1}}{i+1},
\eeq
and\\
\beq\label{eq2}
I-\sum_{i=0}^{2n}\frac{p_i(\theta)(2a)^{\alpha(i+1)}(d^{-\frac{1}{2\beta}}I^{\frac{1}{2\beta}}-R)^{\alpha(i+1)}\{C(\frac{1}{2}(\tau-[\tau])T_0)\}^{i+1}}{i+1}=d(d^{-\frac{1}{2\beta}}I^{\frac{1}{2\beta}}-R)^{2\beta}.
\eeq

Using the Taylor expansion of $f_{\sigma}(x)=(1+x)^{\sigma}=1+\sigma x+x^2\int_0^1f_{\sigma}''(sx)(1-s)ds,$ where $f_{\sigma}''(sx)=\sigma(\sigma-1)(1+sx)^{\sigma-2},$ and with the constants $\alpha,\ \beta,\ d,$ we have
\beq\label{eq3}
d(d^{-\frac{1}{2\beta}}I^{\frac{1}{2\beta}}-R)^{2\beta}=I(1-d^{\frac{1}{2\beta}}I^{-\frac{1}{2\beta}}R)^{\frac{1}{2\beta}}=I-2\beta d^{\frac{1}{2\beta}}I^{1-\frac{1}{2\beta}}R+IR_1,
\eeq
where
\beq\label{eq4} R_1=(d^{\frac{1}{2\beta}}I^{-\frac{1}{2\beta}}R)^2\int_0^1f''_{2\beta}(sd^{\frac{1}{2\beta}}I^{-\frac{1}{2\beta}}R)(1-s)ds,
\eeq
and
\beq\label{eq5}
\begin{array}{ll}
\displaystyle
I-\sum_{i=0}^{2n}\frac{p_i(\theta)(2a)^{\alpha(i+1)}(d^{-\frac{1}{2\beta}}I^{\frac{1}{2\beta}}-R)^{\alpha(i+1)}\{C(\frac{1}{2}(\tau-[\tau])T_0)\}^{i+1}}{i+1}=\\
\displaystyle I-\sum_{i=0}^{2n}\mu_i(\theta,\tau)I^{\frac{\alpha(i+1)}{2\beta}}(1-\alpha(i+1)d^{\frac{1}{2\beta}}I^{-\frac{1}{2\beta}}R+R_2),
\end{array}
\eeq
where
\beq\label{eq6}
\mu_i(\theta,\tau)=\frac{p_i(\theta)(2a)^{\alpha(i+1)}d^{-\frac{\alpha(i+1)}{2\beta}}\{C(\frac{1}{2}(\tau-[\tau])T_0)\}^{i+1}}{i+1},
\eeq
\beq\label{eq7}
R_2=(d^{\frac{1}{2\beta}}I^{-\frac{1}{2\beta}}R)^2\int_0^1f''_{\alpha(i+1)}(sd^{\frac{1}{2\beta}}I^{-\frac{1}{2\beta}}R)(1-s)ds.
\eeq

From (\ref{eq2})-(\ref{eq7}), we have
\beq\label{eq8}
R(I,\theta,\tau)=g_1(I,\theta,\tau)+g_2(I,\theta,\tau)+g_3(I,\theta,\tau)+g_4(I,\theta,\tau),
\eeq
where
\beq\label{eq9}
g_1(I,\theta,\tau)=\sum_{i=0}^{2n}(2\beta d^{\frac{1}{2\beta}})^{-1}\mu_i(\theta,\tau)I^{\frac{\alpha(i+1)}{2\beta}-1},
\eeq
\beq\label{eq10}
g_2(I,\theta,\tau)=-\sum_{i=0}^{2n}(2\beta)^{-1}\mu_i(\theta,\tau)I^{\frac{\alpha(i+1)+1}{2\beta}-1}R,
\eeq
\beq\label{eq11}
g_3(I,\theta,\tau)=\sum_{i=0}^{2n}(2\beta d^{\frac{1}{2\beta}})^{-1}\mu_i(\theta,\tau)I^{\frac{\alpha(i+1)+1}{2\beta}-1}R_2,
\eeq
\beq\label{eq12}
g_4(I,\theta,\tau)=(2\beta d^{\frac{1}{2\beta}})^{-1}I^{\frac{1}{2\beta}}R_1.
\eeq

For $I$ large enough, we can prove the following inequality hold true,
\beq\label{ineq5}
|R|<.I^{\frac{1}{2}},
\eeq
from (\ref{eq8})-(\ref{eq12}).

Assume we have proved (\ref{ineq1}) for $j+k\leq 4$ already, i.e., for I large enough,
\beq\label{ineq6}
|D_{I}^{j}D_{\theta}^{k}R|<.I^{\frac{1}{2}-j},\quad {\rm for}\quad 0\leq j+k\leq4.
\eeq

We are going to prove $|D_{I}^{j}D_{\theta}^{k}R|<.I^{\frac{1}{2}-j},\quad {\rm for}\quad j+k=5.$ From (\ref{eq8}), we have
\beq\label{eq13}
D_{I}^{j}D_{\theta}^{k}R=D_{i}^{j}D_{\theta}^{k}g_1+D_{I}^{j}D_{\theta}^{k}g_2+D_{i}^{j}D_{\theta}^{k}g_3+D_{i}^{j}D_{\theta}^{k}g_4.
\eeq

We will give the estimates of every term of the right hand of Eq.(\ref{ineq6}) as follows.
\begin{Lemma}
\beq\label{ineq7}
|D_{I}^{j}D_{\theta}^{k}g_1|<.I^{\frac{1}{2}-j},\quad {\rm for}\quad j+k=5.
\eeq
\end{Lemma}
\Proof From (\ref{eq9}), we have
\beq\label{eq14}
D_{I}^{j}D_{\theta}^{k}g_1(I,\theta,\tau)=\sum_{i=0}^{2n}(2\beta d^{\frac{1}{2\beta}})^{-1}D_{I}^{j}D_{\theta}^{k}\mu_i(\theta,\tau)I^{\frac{\alpha(i+1)}{2\beta}-1},
\eeq

Using Leibnitz's rule, we can easily prove
\beq\label{ineq8}
|D_{\theta}^{k}\mu_i(\theta,\tau)I^{\frac{\alpha(i+1)}{2\beta}-1}|<.I^{\frac{\alpha(i+1)}{2\beta}-1-j},\quad {\rm for} j+k=5,
\eeq
from (\ref{eq6}) and the definition of $p_i(\tau)$.

Finally, from (\ref{eq13}), we have
$$
|D_{I}^{j}D_{\theta}^{k}g_1(I,\theta,\tau)|=|\sum_{i=0}^{2n}(2\beta d^{\frac{1}{2\beta}})^{-1}D_{I}^{j}D_{\theta}^{k}\mu_i(\theta,\tau)I^{\frac{\alpha(i+1)}{2\beta}-1}|<.I^{\frac{1}{2}-j}.
$$

\begin{Lemma}
\beq\label{ineq9}
|D_{I}^{j}D_{\theta}^{k}g_2|<.I^{-\frac{n+2}{2n+2}-j}+I^{-\frac{1}{2n+2}}|D_{I}^{j}D_{\theta}^{k}R|,\quad {\rm for}\quad j+k=5.
\eeq
\end{Lemma}
\Proof From (\ref{eq10}), we have
\beq\label{eq14}
D_{i}^{j}D_{\theta}^{k}g_2(I,\theta,\tau)=-\sum_{i=0}^{2n}(2\beta)^{-1}D_{i}^{j}D_{\theta}^{k}\mu_i(\theta,\tau)I^{\frac{\alpha(i+1)+1}{2\beta}-1}R,
\eeq

Using Leibnitz's rule, we have
\beq\label{eq15}
\begin{array}{ll}
\displaystyle
D_{I}^{j}D_{\theta}^{k}\mu_i(\theta,\tau)I^{\frac{\alpha(i+1)}{2\beta}-1}R=\sum_{\mbox{\tiny$\begin{array}{c}
j_1+j_2=j,\\
k_1+k_2=k,\\
j+k=5\end{array}$}}D_{I}^{j_1}D_{\theta}^{k_1}(\mu_i(\theta,\tau)I^{\frac{\alpha(i+1)}{2\beta}-1})D_{I}^{j_2}D_{\theta}^{k_2}R\\
\displaystyle =\sum_{\mbox{\tiny$\begin{array}{c}
j_1+j_2=j,\\
k_1+k_2=k,\\
j+k\leq4\end{array}$}}D_{I}^{j_1}D_{\theta}^{k_1}(\mu_i(\theta,\tau)I^{\frac{\alpha(i+1)}{2\beta}-1})D_{I}^{j_2}D_{\theta}^{k_2}R+\mu_i(\theta,\tau)I^{\frac{\alpha(i+1)}{2\beta}-1}D_{I}^{j}D_{\theta}^{k}R.
\end{array}
\eeq

Then we can prove that
$$
|D_{I}^{j}D_{\theta}^{k}g_2|<.I^{-\frac{n+2}{2n+2}-j}+I^{-\frac{1}{2n+2}}|D_{I}^{j}D_{\theta}^{k}R|,\quad {\rm for}\quad j+k=5,
$$
from (\ref{eq6}) and (\ref{ineq6}).

\begin{Lemma}
\beq\label{ineq10}
|D_{I}^{j}D_{\theta}^{k}g_3|<.I^{\frac{n-1}{2n+2}-j}+I^{-\frac{2}{2n+2}}|D_{I}^{j}D_{\theta}^{k}R|,\quad {\rm for}\quad j+k=5.
\eeq
\end{Lemma}

\Proof From (\ref{eq12}), we have
\beq\label{eq16}
D_{I}^{j}D_{\theta}^{k}g_3(I,\theta,\tau)=\sum_{i=0}^{2n}(2\beta d^{\frac{1}{2\beta}})^{-1}D_{I}^{j}D_{\theta}^{k}\mu_i(\theta,\tau)I^{\frac{\alpha(i+1)+1}{2\beta}-1}R_2,
\eeq

Firstly, we give some claims before proving (\ref{ineq10}). Similar to \cite{[WW]} Lemma 3.3, we have
\begin{Lemma} Let $F:\RR\rightarrow\RR$ and $f:\RR^2\rightarrow\RR$ be sufficiently smooth, then for $j+k\geq1$
\beq\label{eq17}
D_{x}^{j}D_{y}^{k}[F\circ f]=\sum_{\mbox{\tiny$\begin{array}{c}
1\leq s\leq j+k,\\
\vec{j}=(j_1,...,j_s),|\vec{j}|=j,\\
\vec{k}=(k_1,...,k_s),|\vec{k}|=k\end{array}$}}C_{s,\vec{j},\vec{k}}(D_{f}^{s}F(f))(D_{x}^{j_1}D_{y}^{k_1}f)...(D_{x}^{j_s}D_{y}^{k_s}f)
\eeq
\end{Lemma}
where $|\vec{j}|=j_1+...+j_s,\ |\vec{k}|=k_1+...+k_s.$

From (\ref{eq7}), we have
\beq\label{eq18}
D_{I}^{j}D_{\theta}^{k}R_2=D_{I}^{j}D_{\theta}^{k}(d^{\frac{1}{2\beta}}I^{-\frac{1}{2\beta}}R)^2\int_0^1f''_{\alpha(i+1)}(sd^{\frac{1}{2\beta}}I^{-\frac{1}{2\beta}}R)(1-s)ds.
\eeq
and
\beq\label{ineq11}
\begin{array}{c}
|D_{I}^{j}D_{\theta}^{k}\int_0^1f''_{\alpha(i+1)}(sd^{\frac{1}{2\beta}}I^{-\frac{1}{2\beta}}R)(1-s)ds|<.|D_{I}^{j}D_{\theta}^{k}f''_{\alpha(i+1)}(sd^{\frac{1}{2\beta}}I^{-\frac{1}{2\beta}}R)|\\ <.|D_{I}^{j}D_{\theta}^{k}(1+sd^{\frac{1}{2\beta}}I^{-\frac{1}{2\beta}}R)^{\sigma-2}|\equiv|D_{I}^{j}D_{\theta}^{k}F_{\sigma}(I^{-\frac{1}{2\beta}}R)|\ \ \ \ \ \ \ \ \ \ \ \ \ \ \ \ \ \ \  \ \ \  \
\end{array}
\eeq

Let $F_{\sigma}(x)=(1+sd^{\frac{1}{2\beta}}x)^{\sigma-2},$ then we have
\beq\label{eq19}
D_{x}^{s}F_{\sigma}(x)=\frac{d^{s}}{dx^{s}}F_{\sigma}(x)=C_{\sigma s}(1+sd^{\frac{1}{2\beta}}x)^{\sigma-2-s}
\eeq
where $C_{\sigma s}$ are constants.

With (\ref{eq17}) and (\ref{eq19}) and  $|I^{-\frac{1}{2\beta}}R|<.\frac{1}{2}$, provide $I$ large enough, we obtain that
\beq\label{ineq12}
\begin{array}{ll}
|D_{I}^{j}D_{\theta}^{k}[F(I^{-\frac{1}{2\beta}}R)]|\\
\displaystyle
=|\sum_{\mbox{\tiny$\begin{array}{c}
1\leq s\leq j+k,\\
\vec{j}=(j_1,...,j_s),|\vec{j}|=j,\\
\vec{k}=(k_1,...,k_s),|\vec{k}|=k\end{array}$}}C_{s,\vec{j},\vec{k}}(D^{s}F(I^{-\frac{1}{2\beta}}R))(D_{I}^{j_1}D_{\theta}^{k_1}(I^{-\frac{1}{2\beta}}R))...(D_{I}^{j_s}D_{\theta}^{k_s}(I^{-\frac{1}{2\beta}}R))|\\
\displaystyle
<.\sum_{\mbox{\tiny$\begin{array}{c}
1\leq s\leq j+k,\\
\vec{j}=(j_1,...,j_s),|\vec{j}|=j,\\
\vec{k}=(k_1,...,k_s),|\vec{k}|=k\end{array}$}}|(D_{I}^{j_1}D_{\theta}^{k_1}(I^{-\frac{1}{2\beta}}R))|...|(D_{I}^{j_s}D_{\theta}^{k_s}(I^{-\frac{1}{2\beta}}R))|
\end{array}
\eeq

For $j+k\leq4$, we have
\beq\label{ineq13}
|D_{I}^{j}D_{\theta}^{k}(I^{-\frac{1}{2\beta}}R)|<.I^{-\frac{1}{2n+2}-j}.
\eeq

When $j+k=5$,
\beq\label{ineq14}
|D_{I}^{j}D_{\theta}^{k}(I^{-\frac{1}{2\beta}}R)|<.I^{-\frac{1}{2n+2}-j}+I^{-\frac{n+2}{2n+2}}|D_{I}^{j}D_{\theta}^{k}R|.
\eeq

From (\ref{eq19})-(\ref{ineq13}), we can prove that
\beq\label{ineq15}
\begin{array}{ll}
|D_{I}^{j}D_{\theta}^{k}[F_{\sigma}(I^{-\frac{1}{2\beta}}R)]|\\
<.I^{-\frac{2}{2n+2}-j}+I^{-\frac{3}{2n+2}-j}+I^{-\frac{4}{2n+2}-j}+I^{-\frac{5}{2n+2}-j}+I^{-\frac{n+2}{2n+2}}|D_{I}^{j}D_{\theta}^{k}R|\\
<.I^{-\frac{2}{2n+2}-j}+I^{-\frac{n+2}{2n+2}}|D_{I}^{j}D_{\theta}^{k}R|\quad {\rm for}\quad j+k=5.
\end{array}
\eeq

Thus, let $\sigma=\alpha(i+1)$, we have
\beq\label{ineq15}
|D_{I}^{j}D_{\theta}^{k}[F_{\alpha(i+1)}(I^{-\frac{1}{2\beta}}R)]|<.I^{-\frac{2}{2n+2}-j}+I^{-\frac{n+2}{2n+2}}|D_{I}^{j}D_{\theta}^{k}R|\quad {\rm for}\quad j+k=5.
\eeq

Using Leibnitz's rule again and with (\ref{ineq13})-(\ref{ineq15}), we can obtain
\beq\label{ineq16}
|D_{I}^{j}D_{\theta}^{k}g_3|<.I^{-\frac{2}{2n+2}-j}+I^{-\frac{n+3}{2n+2}}|D_{I}^{j}D_{\theta}^{k}R|\quad {\rm for}\quad j+k\leq5.
\eeq

Finally, we prove that
\beq\label{ineq17}
|D_{I}^{j}D_{\theta}^{k}g_3(I,\theta,\tau)|=|\sum_{i=0}^{2n}(2\beta d^{\frac{1}{2\beta}})^{-1}D_{I}^{j}D_{\theta}^{k}\mu_i(\theta,\tau)I^{\frac{\alpha(i+1)+1}{2\beta}-1}R_2|<.I^{\frac{n-1}{2n+2}-j}+I^{-\frac{2}{2n+2}}|D_{I}^{j}D_{\theta}^{k}R|,
\eeq
\begin{Lemma}
\beq\label{ineq18}
|D_{I}^{j}D_{\theta}^{k}g_4|<.I^{\frac{n}{2n+2}-j}+I^{-\frac{1}{2n+2}}|D_{I}^{j}D_{\theta}^{k}R|\quad {\rm for}\quad j+k=5.
\eeq
\end{Lemma}
The proof is similar to Lemma 4.3, we omit it.

Now, with lemma 4.1, 4.2, 4.3, 4.5 and (\ref{ineq6}), we have
\beq\label{semi}
\begin{array}{ll}
|D_{I}^{j}D_{\theta}^{k}R|\\
<.I^{\frac{1}{2}-j}+I^{-\frac{n+2}{2n+2}-j}+I^{-\frac{1}{2n+2}}|D_{I}^{j}D_{\theta}^{k}R|+I^{\frac{n-1}{2n+2}-j}+I^{-\frac{2}{2n+2}}|D_{I}^{j}D_{\theta}^{k}R|+I^{\frac{n}{2n+2}-j}+I^{-\frac{1}{2n+2}}|D_{I}^{j}D_{\theta}^{k}R|\\
<.I^{\frac{1}{2}-j}\quad {\rm for}\quad j+k=5.
\end{array}
\eeq

Thus, we have $|D_{I}^{j}D_{\theta}^{k}R|<.I^{\frac{1}{2}-j}\quad {\rm provided}\quad 0\leq j+k\leq5.$ The proof of Lemma 3.1 is completed.

\vskip 2mm
\section{Conclusion}

\ \ \ \ \ In this paper, we study the Lagrangian stability for a class of impact oscillators with time dependent polynomial potentials. By exchanging the roles of angle and time, we overcome the nonsmoothness on angle variable caused by the existence of impact. Then we obtain an equivalent Hamiltonian which is sufficiently smooth in new angle variable by some canonical transformations. Thus by direct application of Moser's twist theorem, we obtain the existence of invariant curves for the Poincar\'{e} mapping of the equations, which implies the Lagrangian stability.

\end{document}